\documentclass[12pt]{article}
\usepackage{amssymb,amsmath,amsopn,amsfonts}
\usepackage[dvips]{color}
\usepackage{colordvi,multicol}
\usepackage{epsf}
\newfam\msbfam
\font\tenmsb=msbm10 \textfont\msbfam=\tenmsb \font\sevenmsb=msbm7
\scriptfont\msbfam=\sevenmsb \font\fivemsb=msbm5
\scriptscriptfont\msbfam=\fivemsb

\def\th#1{\vspace{1mm}\noindent{\bf #1}\quad}

\def\proof{\vspace{1mm}\noindent{\it Proof}\quad}

\def\bc{\begin{center}}
\def\ec{\end{center}}
\def\no{\noindent}
\def\hang{\hangindent\parindent}
\def\textindent#1{\indent\llap{\qquad #1\ \ \enspace}\ignorespaces}
\def\ref{\par\hang\textindent}

\begin{document}

\title{ {\bf BV functions in a Gelfand triple
 for differentiable measure and its applications
\thanks{Research supported  the DFG through IRTG 1132 and CRC 701}\\} }
\author{{\bf Michael R\"{o}ckner}$^{\mbox{c},}$, {\bf Rongchan Zhu}$^{\mbox{b,c},}$, {\bf Xiangchan Zhu}$^{\mbox{a,c},}$
\date {}
\thanks{E-mail address:  roeckner@math.uni-bielefeld.de(M. R\"{o}ckner), zhurongchan@126.com(R. C. Zhu), zhuxiangchan@126.com(X. C. Zhu)}\\ \\
$^{\mbox{a}}$School of Sciences, Beijing Jiaotong University, Beijing 100044, China\\
$^{\mbox{b}}$Department of Mathematics, Beijing Institute of Technology, Beijing 100081,
 China,\\
$^{\mbox{c}}$ Department of Mathematics, University of Bielefeld, D-33615 Bielefeld, Germany,}

\maketitle

\begin{abstract}
\vskip 0.1cm \noindent In this paper, we introduce a definition of BV functions for (non-Gaussian) differentiable measure in a Gelfand triple which is an extension of the definition of BV functions in [RZZ12], using Dirichlet form theory. By this definition, we can analyze the reflected stochastic  quantization problem associated with a self-adjoint operator $A$ and a cylindrical Wiener process on a convex set $\Gamma$ in a Banach space $E$. We prove the existence of a martingale solution of this problem if $\Gamma$ is a regular convex set.
\end{abstract}
\vspace{1mm}
\no{\footnotesize{\bf 2000 Mathematics Subject Classification AMS}:\hspace{2mm} 60Gxx, 31C25, 60G60, 26A45}
 \vspace{2mm}

\no{\footnotesize{\bf Keywords}:\hspace{2mm}  Dirichlet
forms, Stochastic reflection problems, BV function, Gelfand triples, Integration by parts formula in infinite dimensions, Differentiable measure, Stochastic quantization}

\vspace{1mm}
\section{Introduction}
 The theory of bounded variation (BV) functions in infinite dimensions was developed a lot in the last years (cf. [Fu00], [FH01], [AMMP10], [Hi10], [ADP10], [ADGP12], [RZZ12]).
 A definition of BV functions in abstract Wiener spaces has been given by M. Fukushima in [Fu00], and M. Fukushima and M. Hino in [FH01], based upon Dirichlet form theory. In [RZZ12] the authors introduce BV functions in a Gelfand triple, which is an extension of BV functions in a Hilbert space defined as in [ADP10]. In [ADGP12], the authors define BV functions in a Hilbert space with respect to log-concave measure, which are absolutely continuous with respect to a Gaussian measure, and give a characterization in terms of the semigroup of a SDE.

 In this paper, realizing that to have integration by parts in sufficiently many directions is sufficient for the definition of BV functions, we analogously define BV functions on a Banach space $E$ and replace the Gaussian measure with a differentiable measure $\mu$ in a Gelfand triple (see Definitions 2.1 and 3.1 below). Differentiable measures (see [Bo10]) form a general class which contains besides Gaussian measures also measures which are not absolutely continuous with respect to any Gaussian measures, such as Gibbs measure from statistical mechanics and Euclidean field theory (see e.g. [AKR97a], [AKR97b]). This definition of BV functions can be seen as an extension of BV functions in a Gelfand triple and BV functions in abstract Wiener spaces, but since we use differentiable measure, as a price to pay, we need the BV functions to belong to $L^2(E,\mu)$. Here we use a version of the Riesz-Markov representation theorem in infinite dimensions proved by M. Fukushima using the quasi-regularity of the Dirichlet form (see [MR92]) to give a characterization of BV functions.

 We consider the Dirichlet form $$\mathcal{E}^\rho(u,v)=\frac{1}{2}\sum_{k=1}^\infty\int_E\frac{\partial u}{\partial e_k}\frac{\partial v}{\partial e_k}\rho d\mu,$$
(where $E$ is a Banach space with a Hilbert space $H\subset E$ continuously and densely, $e_j,j\in \mathbb{N}$ is an orthonormal basis in $H$, $\mu$ is a differentiable measure in $E$ and $\rho$ is a BV function)
and its associated Markov process. Using BV functions, we obtain a Skorohod-type representation for the associated process, if $\rho=I_\Gamma$ and $\Gamma$ is a convex set (see Theorem 3.4 below).

As a consequence of these results, we can solve the following stochastic differential inclusion in the Banach space $E$:
$$\left\{\begin{array}{ll}dX(t)+(AX(t)+:p(X):+N_\Gamma(X(t)))dt\ni dW(t),&\ \ \ \ \textrm{ }\\X(0)=x,&\ \ \ \ \textrm{ } \end{array}\right.\eqno(1.1)$$
 if $\Gamma$ is regular.
Here $A:D(A)\subset H\rightarrow H$ is a self-adjoint operator. $N_\Gamma(x)$ is the normal cone to $\Gamma$ at $x$ and $W(t)$ is a cylindrical Wiener process in $H$. The solution to (1.1) is called  reflected stochastic quantization process. We would like to stress that our results apply to models from Euclidean 2D-quantum fields ("$P(\phi)_2$-models") both in finite (cf. Theorems 4.1.1-4.1.4) and infinite volume (cf. Theorem 4.2.1 and Theorem 4.2.2). The latter is generally much more difficult than the first.

This kind of reflection problem without the term involving $:p(X):$, i.e. $p=0$, in infinite dimension was first studied (strongly solved) in [NP92], when $H=L^2(0,1)$, $A$ is the Laplace operator with Dirichlet or Neumann boundary conditions and $\Gamma$ is the convex set of all nonnegative functions of $L^2(0,1)$; see also [Za02]. In [BDT09] the authors study the situation when $p=0$ and $\Gamma$ is a regular convex set  with nonempty interior. They get precise information about the corresponding Kolmogorov operator. In [RZZ12], by using BV functions we deduce that, if $p=0$, (1.1) has a unique strong solution in the probabilistic sense. In the above references, the reflected problem was considered only with space dimension 1. By using the BV functions, we obtain the martingale solutions to reflected OU process with space dimension 2.

Following the pioneering paper of [JM85] in finite volume, the stochastic quantization problem with space dimension 2 (without reflection term) was studied in [AR89, 90] ("infinite and finite volume"), [AR91]("infinite and finite volume"), [RZ92]("finite volume"), [LR98]("finite volume") by using Dirichlet form theory (see also [DT00] for an approach via SPDE, also in finite volume). Da Prato and Debussche in [DD03] proved the existence and uniqueness of a strong solution of this problem, but only in the finite volume case. By using BV functions, in the present paper we obtain martingale solutions to the reflected stochastic quantization problem in finite and infinite volume.

This paper is organized as follows. In Section 2, we consider the appropriate Dirichlet form
and its associated "distorted" process. We introduce BV functions for differentiable measures in Section 3, and prove the Skorohod type representation for the distorted process. In Section 4, we give examples of BV functions for differentiable measure, in particular, the reflected stochastic  quantization problem mentioned above.

\section{The Dirichlet form and the associated distorted process}
 Let $E$ be a separable Banach space and $H$ be a real separable Hilbert space (with scalar product $\langle\cdot,\cdot\rangle$ and norm denoted by $|\cdot|$), continuously and densely embedded in $E$. We denote their Borel $\sigma$-algebras by $\mathcal{B}(H), \mathcal{B}(E)$ respectively.
 Here identifying $H$ with its dual $H^*$ we obtain the continuous and dense embeddings
$$E^*\subset H(\equiv H^*)\subset E.$$
 It follows that
 $$_{E^*}\langle z,v\rangle_E=\langle z,v\rangle_H \forall z\in E^*,v\in H.$$

For a (positive) measure $\mu$ on $\mathcal{B}(E)$ let $L^p(E,\mu), p\in [1,\infty]$, denote the corresponding real $L^p$-spaces equipped with the usual norms $\|\cdot\|_p$.

\vskip.10in

Let $$\mathcal{F}C_b^1=\{u:u(z)=f({ }_{E^*}\!\langle l_1,z\rangle_E,{ }_{E^*}\!\langle l_2,z\rangle_E,...,{ }_{E^*}\!\langle l_m,z\rangle_E),z\in E, l_1,l_2,...,l_m\in E^*, m\in \mathbb{N}, f\in C_b^1(\mathbb{R}^m)\}.$$
 Define for $u\in \mathcal{F}C_b^1$ and $l\in H$, $$\frac{\partial u}{\partial l}(z):=\frac{d}{ds}u(z+sl)|_{s=0},z\in E,$$
  that is, by the chain rule,
  $$\frac{\partial u}{\partial l}(z)=\sum_{j=1}^m\partial_jf({ }_{E^*}\!\langle l_1,z\rangle_E,{ }_{E^*}\!\langle l_2,z\rangle_E,...,{ }_{E^*}\!\langle l_m,z\rangle_E)\langle l_j,l\rangle.$$
Let $Du$ denote the $H$-derivative of $u\in \mathcal{F}C_b^1$, i.e. the map from $E$ to $H$ such that $$\langle Du(z),l\rangle=\frac{\partial u}{\partial l}(z)\textrm{ for all } l\in H, z\in E$$
Let $\mu(\neq0)$ be a finite positive Radon measure on $\mathcal{B}(E)$ having the following property: if a function $\varphi\in \mathcal{F}C_b^1$ is equal to zero $\mu$-almost everywhere, then $\frac{\partial\varphi}{\partial l}=0$ $\mu$-almost everywhere for all $l\in H$. In particular, this holds for a measure with full topological support, i.e. $\mu(U)>0$ for all nonempty open $U\subset E$. Now we recall the following definition from [Fo66, Fo68].
\vskip.10in
\th{Definition 2.1} A measure $\mu$ on $E$ is called differentiable along $v\in E$ in the sense of Fomin if there exists a finite signed measure $d_v\mu$ such that for any $\varphi \in\mathcal{F}C_b^1$ the following equality holds:
$$\int\frac{\partial\varphi(x)}{\partial v}\mu(dx)=-\int \varphi(x)d_v\mu(dx),$$
and $d_v\mu$ is absolutely continuous with respect to $\mu$. Its density
 with respect to $\mu$ will be denoted by $\beta_v$.

\th{Remark} By Corollary 3.3.2, Proposition 3.3.4 and Corollary 3.3.5 in  [Bo10], this is equivalent to Definition 3.1.1 in [Bo10]: A measure $\mu$ on $E$ is called differentiable along $v\in E$ in the sense of Fomin if, $d_v\mu$ is a signed measure of bounded variation such that for every set $A\in\mathcal{B}(E) $, there exists the finite limit
$$d_v\mu(A):=\lim_{t\rightarrow0}\frac{\mu(A+tv)-\mu(A)}{t}.$$

Set
$$H(\mu):=\{h\in E: \mu \textrm{ is differentiable along }h\textrm{ and }\|\beta_h\|_2<\infty \},$$
endowed with the norm $\|h\|_{H(\mu)}:=\|\beta_h\|_2$. Then $H(\mu)$ is a Hilbert space continuously embedded in E.

From now on we fix a differentiable measure $\mu$ on $E$ satisfying the following:

\th{Hypothesis 2.2} Let $Q: H\rightarrow H$ be a nonnegative definite  linear bounded symmetric operator such that $Q^{1/2}(H)\subset H(\mu)$ and there exists an orthonormal basis $\{e_j\}$ of $H$ consisting of eigen-functions for $Q$ with corresponding eigenvalues $\lambda_j\in (0,\infty),j\in\mathbb{N},$ that is,
$$Qe_j=\lambda_je_j,j\in\mathbb{N},$$ and such that $\{e_j\}\subset E^*$.

For $\rho\in L^1_+(E,\mu)$ we consider
$$\mathcal{E}^\rho(u,v)=\frac{1}{2}\sum_{k=1}^\infty\int_E\frac{\partial u}{\partial e_k}\frac{\partial v}{\partial e_k}\rho d\mu,u,v\in \mathcal{F}C_b^1,$$
where $ L^1_+(E,\mu)$ denotes the set of all non-negative elements in $ L^1(E,\mu)$. Let $QR(E)$ be the set of all functions $\rho\in L^1_+(E,\mu)$ such that $(\mathcal{E}^\rho, \mathcal{F}C_b^1)$ is closable on $L^2(F,\rho\cdotp\mu$), where $F:=Supp[\rho\cdot\mu]$. Its closure is denoted by $(\mathcal{E}^\rho, \mathcal{F}^\rho)$. We denote by $\mathcal{F}_e^\rho$ the extended Dirichlet space of $(\mathcal{E}^\rho, \mathcal{F}^\rho)$, that is,  $u\in \mathcal{F}_e^\rho$ if and only if $|u|<\infty \textrm{ }\rho\cdot\mu-a.e.$ and there exists a sequence $\{u_n\}$ in $\mathcal{F^\rho}$ such that $\mathcal{E}^\rho(u_m-u_n,u_m-u_n)\rightarrow0$ as $n\geq m\rightarrow\infty$ and $u_n\rightarrow u \textrm{  } \textrm{  }\rho\cdot\mu-a.e.$ as $n\rightarrow \infty$.
\vskip.10in
\th{Theorem 2.3} Let $\rho\in QR(E)$. Then $(\mathcal{E}^\rho, \mathcal{F}^\rho)$ is  a quasi-regular local Dirichlet form on $L^2(F;\rho\cdot\mu)$ in the sense of [MR92, IV Definition 3.1].

\proof The assertion follows from the main result in [RS92]. $\hfill\Box$
\vskip.10in

By virtue of Theorem 2.3 and [MR92], there exists a diffusion process $M^\rho=(\Omega,\mathcal{M},\{\mathcal{M}_t\},\theta_t,$ $X_t,P_z)$ on $F$ associated with the Dirichlet form $(\mathcal{E}^\rho, \mathcal{F}^\rho).$ $M^\rho$ will be called distorted process on $F$. Since constant functions are in $\mathcal{F}^\rho$ and $\mathcal{E}^\rho(1,1)=0$, $M^\rho$ is recurrent and conservative. By $\textbf{A}_+^\rho$ we denote the set of all positive continuous additive functionals (PCAF in abbreviation) of $M^\rho$, and define $\textbf{A}^\rho:=\textbf{A}^\rho_+-\textbf{A}^\rho_+$. For $A\in \textbf{A}^\rho$, its total variation process is denoted by $\{A\}$. We also define $\textbf{A}^\rho_0:=\{A\in \textbf{A}^\rho|E_{\rho\cdotp\mu}(\{A\}_t)<\infty, \forall t>0\}$. Each element in $\textbf{A}^\rho_+$ has a corresponding positive $\mathcal{E}^\rho$-smooth measure on $F$ by the Revuz correspondence. The set of all such measures will be denoted by $S^\rho_+$. Accordingly, $A_t\in\textbf{A}^\rho$ corresponds to a $\nu\in S^\rho:=S^\rho_+-S^\rho_+$, the set of all $\mathcal{E}^\rho$-smooth signed measure in the sense that   $A_t=A_t^1-A_t^2$ for $A_t^k\in \textbf{A}^\rho_+,k=1,2$ whose  Revuz measures are $\nu^k, k=1,2$ and $\nu=\nu^1-\nu^2$ is the Hahn-Jordan decomposition of $\nu$ . The element of $\textbf{A}^\rho$ corresponding to $\nu\in S^\rho$ will be denoted by$A^\nu$.

Note that for each $l\in E^*$ the function $u(z)={ }_{E^*}\!\langle l,z\rangle_E$ belongs to the extended Dirichlet space $\mathcal{F}^\rho_e$ and $$\mathcal{E}^\rho(u(\cdot),v)=\frac{1}{2}\int\frac{\partial v(z)}{\partial l}\rho(z)d\mu(z)\textrm{ }\forall v\in \mathcal{F}C_b^1.\eqno(2.1)$$
On the other hand, the AF $_{E^*}\langle l,X_t-X_0\rangle_E$ of $M^\rho$ admits a unique decomposition into a sum of a martingale AF ($M_t$) of finite energy and a CAF ($N_t$) of zero energy. More precisely, for every $l\in E^*$, $$_{E^*}\langle l,X_t-X_0\rangle_E=M^l_t+N^l_t\textrm{ }\forall t\geq0\textrm{ } P_z- a.s.\eqno(2.2)$$ for $\mathcal{E}^\rho$-q.e. $z\in F$.

Now for $\rho\in L^1(E,\mu)$ and $l\in E^*$, we say that $\rho\in BV_l(E)$ if there exists a constant $C_l>0$,
$$|\int_E\frac{\partial v(z)}{\partial l}\rho(z)d\mu(z)|\leq C_l\parallel v\parallel_\infty\textrm{ }\forall v\in \mathcal{F}C_b^1.\eqno(2.3)$$

By the same argument as in [FH01, Theorem 2.1], we obtain the following:
\vskip.10in
\th{Theorem 2.4} Let $\rho\in L^1_+$ and $l\in E^*$.

(1) The following two conditions are equivalent:

(i)$\rho\in BV_l(E)$

(ii) There exists a (unique) signed measure $\nu_l$ on $F$ of finite total variation such that
$$\frac{1}{2}\int\frac{\partial v(z)}{\partial l}\rho(z)d\mu(z)=-\int_Fv(z)\nu_l(dz) \textrm{ }\forall v\in \mathcal{F}C_b^1.\eqno(2.4)$$
In this case, $\nu_l$ necessarily belongs to $S^{\rho+1}$.

Suppose further that $\rho\in QR(E)$. Then the following condition is also equivalent to the above:

(iii)$N^l\in \textbf{A}_0^\rho$

In this case, $\nu_l\in S^\rho$, and $N^l=A^{\nu_l}$

(2) $M^l$ is a martingale AF with quadratic variation process
$$\langle M^l\rangle_t=t|l|^2,t\geq0.\eqno(2.5)$$
\vskip.10in
\th{Remark 2.5} Recall that the Riesz representation theorem, saying that every positive linear functional on continuous functions can be represented by a measure, is not applicable to obtain Theorem 2.4, $(i)\Rightarrow(ii)$, because of the lack of local compactness of $E$. However, the quasi-regularity of the Dirichlet form provides a means to circumvent this difficulty.

\section{BV functions and distorted processes in $F$}

Let $c_j ,j\in \mathbb{N}$, be a sequence in $[1,\infty)$. Define for $e_j,j\in \mathbb{N}$, as in Hypothesis 2.2
 $$H_1:=\{x\in H|\sum_{j=1}^\infty \langle x,e_j\rangle ^2c_j^2<\infty\},$$
equipped with the inner product
$$\langle x,y\rangle_{H_1}:=\sum_{j=1}^\infty c_j^2\langle x,e_j\rangle \langle y,e_j\rangle.$$
Then clearly $(H_1, \langle,\rangle_{H_1})$ is a Hilbert space such that $H_1\subset H$ continuously and densely. Identifying $H$ with its dual we obtain the continuous and dense embeddings
$$H_1\subset H(\equiv H^*)\subset H^*_1.$$
 It follows that
 $$_{H_1}\langle z,v\rangle_{H_1^*}=\langle z,v\rangle_H \textrm{ }\forall z\in H_1,v\in H,$$ and that
$(H_1,H,H_1^*)$ is a Gelfand triple. Furthermore, $\{\frac{e_j}{c_j}\}$ and $\{c_je_j\}$ are orthonormal bases of $H_1$ and $H_1^*$, respectively. In particular, $e_j\in Q^{1/2}(H)\cap H_1 \cap E^*$.

We also introduce a family of $H$-valued functions on $E$ by
$$(\mathcal{F}C_b^1)_{Q^{1/2}(H)\cap H_1}:=\{G:G(z)=\sum_{j=1}^m g_j(z)l^j,z\in E, g_j\in \mathcal{F}C_b^1,l^j\in Q^{1/2}(H)\cap H_1\}$$
Denote by $D^*$ the adjoint of $D:\mathcal{F}C_b^1\subset L^2(E,\mu)\rightarrow L^2(E,\mu;H)$. That is
$$Dom(D^*):=\{G\in L^2(E,\mu;H)|\mathcal{F}C_b^1\ni u\mapsto \int_E\langle G,Du\rangle d\mu \textrm{ is continuous with respect to } L^2(E,\mu)\}.$$ Obviously, $(\mathcal{F}C_b^1)_{Q^{1/2}(H)\cap H_1}\subset Dom(D^*)$. Then
$$\int_E D^*G(z)f(z)\mu(dz)=\int_E\langle G(z),Df(z)\rangle\mu(dz)\textrm{ } \forall G\in (\mathcal{F}C_b^1)_{Q^{1/2}(H)\cap H_1},f\in \mathcal{F}C_b^1.\eqno(3.1)$$

For $\rho\in L^2(E,\mu)$, we set
 $$V(\rho):=\sup_{G\in (\mathcal{F}C_b^1)_{Q^{1/2}(H)\cap H_1},\|G\|_{H_1}\leq1}\int_E D^*G(z)\rho(z)\mu(dz).\eqno(3.2)$$

\th{Definition 3.1} A function $\rho$ on $E$ is called a BV function in the Gelfand triple $(H_1, H, H^*_1)$($\rho \in BV(H, H_1)$ in notation), if $\rho\in L^2(E,\mu)$ and $V(\rho)$ is finite.
\vskip.10in
We can prove the following theorem by a modification of the proof of [RZZ12, Theorem 3.1].

\th{Theorem 3.2} (i) $BV(H, H_1)\subset\bigcap_{l\in Q^{1/2}(H)\cap H_1\cap E^*} BV_l(E)$.

(ii) Suppose $\rho\in BV(H, H_1)\cap L^1_+(E,\mu)$, then there exist a positive finite measure $\|d\rho\|$ on $E$ and a Borel-measurable map $\sigma_\rho:E\rightarrow H_1^*$ such that $\|\sigma_\rho(z)\|_{H_1^*}=1\textrm{ }\|d\rho\|-a.e, \|d\rho\|(E)=V(\rho)$,
$$\int_ED^*G(z)\rho(z)\mu(dz)=\int_E { }_{H_1}\!\langle G(z),\sigma_\rho(z)\rangle_{H_1^*}\|d\rho\|(dz)\textrm{ }\forall G\in (\mathcal{F}C_b^1)_{Q^{1/2}(H)\cap H_1} \eqno(3.3)$$
and $\|d\rho\|\in S^{\rho+1}$.

Furthermore, if $\rho \in QR(E)$, $\|d\rho\|$ is $\mathcal{E}^\rho$-smooth in the sense that it charges no set of zero $\mathcal{E}_1^\rho$-capacity. In particular, the domain of integration $E$ on both sides of (3.3) can be replaced by $F$, the topological support of $\rho\mu$.

Also, $\sigma_\rho$ and $\|d\rho\|$ are uniquely determined, that is, if there are $\sigma_\rho'$ and $\|d\rho\|'$ satisfying relation (3.3), then $\|d\rho\|=\|d\rho\|'$ and $\sigma_\rho(z)=\sigma_\rho'(z)$ for $\|d\rho\|-a.e.z$

(iii) Conversely, if Eq.(3.3) holds for $\rho\in L^2(E,\mu)$ and for some positive finite measure $\|d\rho\|$ and a map $\sigma_\rho$ with the stated properties, then $\rho\in BV(H,H_1)$ and $V(\rho)=\|d\rho\|(E)$.

(iv) Let $W^{1,2}(E)$ be the domain of the closure of $(D, \mathcal{F}C_b^1)$ with norm
$$\|f\|^2:=\int_E (|f(z)|^2+|Df(z)|^2)\mu(dz).$$
Then $W^{1,2}(E)\subset BV(H,H)$ and  Eq.(3.3) is satisfied for each $\rho\in W^{1,2}(E)$. Furthermore,
$$\|d\rho\|=|D\rho|\cdot\mu, V(\rho)=\int_E|D\rho|\mu(dz),\sigma_\rho=\frac{1}{|D\rho|}D\rho I_{\{|D\rho|>0\}}.$$

\proof (i) Let $\rho\in BV(H,H_1)$. Take $G\in (\mathcal{F}C_b^1)_{Q^{1/2}(H)\cap H_1}$ of the type
$$G(z)=g(z)l,z\in E, g\in \mathcal{F}C_b^1,l\in Q^{1/2}(H)\cap H_1.\eqno(3.4)$$
 By (3.1) and [Bo10, Proposition 3.3.12] $$\aligned\int_ED^*G(z)f(z)\mu(dz)=&\int_E\langle G(z),Df(z)\rangle\mu(dz)\\=&-\int_E\langle l,Dg(z)\rangle f(z)\mu(dz)-\int_E\beta_l(z)g(z)f(z)\mu(dz)\textrm{ }\forall f\in \mathcal{F}C_b^1;\endaligned$$
consequently, $$D^*G(z)=-\langle l, Dg(z)\rangle-g(z)\beta_l(z).\eqno(3.5)$$
 Accordingly,
 $$\int_E\langle l,Dg(z)\rangle \rho(z)\mu(dz)=-\int_ED^*G(z)\rho(z)\mu(dz)-\int_E\beta_l(z) g(z)\rho(z)\mu(dz).\eqno(3.6)$$
 For any $g\in \mathcal{F}C_b^1$, satisfying $\|g\|_\infty\leq1$, by (3.2) the right hand side of (3.6) is dominated by
 $$V(\rho)\|l\|_{H_1}+\|\rho\|_2\|\beta_l\|_2<\infty , $$
 hence, $\rho\in BV_l(H)$.

 (ii) Suppose $\rho\in L^1_+(E,\mu)\bigcap BV(H,H_1)$. By (i) and Theorem 2.4  for each $l\in Q^{1/2}(H)\cap H_1\cap E^*$, there exists a finite signed measure $\nu_l$ on $E$ for which Eq.(2.4) holds. Define $$D^A_l\rho(dz):=2\nu_l(dz)-\beta_l(z)\rho(z)\mu(dz).$$
 In view of (3.6), for any $G$ of  type (3.4), we have
 $$\int_ED^*G(z)\rho(z)\mu(dz)=\int_Eg(z)D_l^A\rho(dz),\eqno(3.7)$$
 which in turn implies $$V(D^A_l\rho)(E)=\sup_{g\in \mathcal{F}C_b^1, \|g\|_\infty\leq1}\int_E g(z)D_l^A\rho(dz)\leq V(\rho)\|l\|_{H_1},\eqno(3.8)$$
 where $V(D^A_l\rho)$ denotes the total variation measure of the signed measure $D_l^A\rho$.

 For the orthonormal basis $\{\frac{e_j}{c_j}\}$ of $H_1$, we set
 $$\gamma_\rho^A:=\Sigma_{j=1}^\infty 2^{-j}V(D_{\frac{e_j}{c_j}}^A\rho),\textrm{ }v_j(z):=\frac{d D_{\frac{e_j}{c_j}}^A\rho(z)}{d\gamma_\rho^A(z)},z\in E, j\in \mathbb{N}.\eqno(3.9)$$
 $\gamma_\rho^A$ is a positive finite measure with $\gamma_\rho^A(E)\leq V(\rho)$  and $v_j$ is Borel-measurable. Since $D_{\frac{e_j}{c_j}}^A\rho$ belongs to $S^{\rho+1}$, so does $\gamma_\rho^A$ . Then for
 $$G_n:=\sum_{j=1}^n g_j\frac{e_j}{c_j}\in(\mathcal{F}C_b^1)_{Q^{1/2}(H)\cap H_1},n\in \mathbb{N},\eqno(3.10)$$
 by (3.7) the following equation holds
 $$\int_ED^*G_n(z)\rho(z)\mu(dz)=\sum_{j=1}^n\int_E g_j(z)v_j(z)\gamma_\rho^A(dz).\eqno(3.11)$$
 Since $|v_j(z)|\leq2^j \textrm{ } \gamma_\rho^A$-a.e. and $\mathcal{F}C^1_b$ is dense in $L^1(E,\gamma_\rho^A)$, we can find $v_{j,m}\in \mathcal{F}C_b^1$ such that
 $$\lim_{m\rightarrow\infty}v_{j,m}=v_j  \textrm{ } \gamma_\rho^A-a.e., $$
 Substituting
 $$g_{j,m}(z):=\frac{v_{j,m}(z)}{\sqrt{\sum_{k=1}^nv_{k,m}(z)^2+1/m}},\eqno(3.12)$$
 for $g_j(z)$ in (3.10) and (3.11) we get
 $$\sum_{j=1}^n\int_E g_{j,m}(z)v_j(z)\gamma_\rho^A(dz)\leq V(\rho),$$
 because $\|G_n(z)\|_{H_1}^2=\sum_{j=1}^n g_{j,m}(z)^2\leq1 \textrm{ }\forall z\in E$. By letting $m\rightarrow\infty$, we obtain $$\int_E\sqrt{\sum_{j=1}^nv_j(z)^2}\gamma_\rho^A(dz)\leq V(\rho)\textrm{ }\forall n\in \mathbb{N}.$$
 Now we define
 $$\| d\rho\|:=\sqrt{\sum_{j=1}^\infty v_j(z)^2}\gamma_\rho^A(dz)\eqno(3.13)$$ and $\sigma_\rho:E\rightarrow H_1^*$ by
 $$\sigma_\rho(z)=\left\{\begin{array}{ll}\sum_{j=1}^\infty\frac{v_j(z)}{\sqrt{\sum_{k=1}^\infty v_k(z)^2}}\cdot c_je_j,&\ \ \ \ \textrm{ if }z\in\{\sum_{k=1}^\infty v_k(z)^2>0\}\\0&\ \ \ \ \textrm{ otherwise. } \end{array}\right.\eqno(3.14)$$
 Then$$\|d\rho\|(E)\leq V(\rho),\textrm{ } \|\sigma_\rho(z)\|_{H_1^*}=1 \textrm{ }\|d\rho\|-a.e.,\eqno(3.15)$$
 $\|d\rho\|$ is $S^{\rho+1}$-smooth and $\sigma_\rho$ is Borel-measurable. By (3.11) we see that the desired equation (3.3) holds for $G=G_n$ as in (3.10).
 It remains to prove (3.3) for any $G$ of type (3.4), i.e. $G=g\cdot l,g\in \mathcal{F}C_b^1,l\in Q^{1/2}(H)\cap H_1$ . In view of (3.6), Eq.(3.3) then reads
  $$-\int_E\langle l, Dg(z)\rangle\rho(z)\mu(dz)-\int_Eg(z)\beta_l(z) \rho(z)\mu(dz)=\int_Eg(z)_{H_1}\langle l,\sigma_\rho(z)\rangle_{H_1^*}\|d\rho\|(dz).\eqno(3.16)$$
 We set $$k_n:=\sum_{j=1}^n\langle l,e_j\rangle e_j=\sum_{j=1}^n\langle l,\frac{e_j}{c_j}\rangle_{H_1} \frac{e_j}{c_j}=\sum_{j=1}^n\langle l,\lambda_j^{1/2}e_j\rangle_{Q^{1/2}(H)} \lambda_j^{1/2}e_j, G_n(z):=g(z)k_n.$$
 Thus $k_n\rightarrow l$ in $H_1$ and $k_n\rightarrow l$ in $Q^{1/2}(H)$ as $n\rightarrow \infty$. Hence by Hypothesis 2.2 and [Bo10, Proposition 5.1.7] $\|\beta_{k_n}-\beta_l\|_2\rightarrow0$.
 But then also
 $$\lim_{n\rightarrow\infty}\int_E \langle Dg,k_n\rangle\rho d\mu=\int_E \langle Dg,l\rangle\rho d\mu,$$
 and
 $$|\int_Eg(z)\beta_{k_n}(z) \rho(z)\mu(dz)-\int_Eg(z)\beta_l(z) \rho(z)\mu(dz)|$$
 $$\leq \|g\|_\infty\|\rho\|_2\|\beta_{k_n}-\beta_l\|_2.$$

 Furthermore, $$\lim_{n\rightarrow\infty}\int_Eg(z)_{H_1}\langle k_n,\sigma_\rho(z)\rangle_{H_1^*}\|d\rho\|(dz)=\int_Eg(z)_{H_1}\langle l,\sigma_\rho(z)\rangle_{H_1^*}\|d\rho\|(dz).$$
 So letting $n\rightarrow\infty$ yields (3.16).

 If $\rho \in QR(E)$, we can get the claimed result by the same arguments as above.

 Uniqueness follows by the same argument as [FH01, Theorem 3.9].

 (iii) Suppose $\rho\in L^2(E,\mu)$ and that Eq.(3.3) holds for some positive finite measure $\|d\rho\|$ and some map $\sigma_\rho$ with the properties stated in (ii). Then clearly
 $$V(\rho)\leq\|d\rho\|(E)$$ and hence $\rho\in BV(H,H_1)$. To obtain the converse inequality, set
 $$\sigma_j(z):=\langle c_je_j,\sigma_\rho(z)\rangle_{H_1^*}=_{H_1}\langle \frac{e_j}{c_j},\sigma_\rho(z)\rangle_{H_1^*}, j\in \mathbb{N}.$$
 Fix an arbitrary $n$. As in the proof of (ii) we can find functions
 $$v_{j,m}\in \mathcal{F}C_b^1, \qquad\lim_{m\rightarrow\infty}v_{j,m}(z)=\sigma_j(z)\textrm{ } \|d\rho\|-a.e.$$
 Define $g_{j,m}(z)$ by (3.12). Substituting $G_{n,m}(z):=\sum_{j=1}^ng_{j,m}(z)\frac{e_j}{c_j}$ for $G(z)$ in (3.3) then yields
 $$\sum_{j=1}^n\int_E g_{j,m}(z)\sigma_j(z)\|d\rho\|(dz)\leq V(\rho).$$
 By letting $m\rightarrow\infty$, we get $$\int_E\sqrt{\sum_{j=1}^n\sigma_j(z)^2}\|d\rho\|(dz)\leq V(\rho)\textrm{ } \forall n\in \mathbb{N}.$$
 We finally let $n\rightarrow\infty$ to obtain $\|d\rho\|(E)\leq V(\rho)$.

 (iv) Obviously the duality relation (3.1) extends to $\rho\in W^{1,2}(E)$ replacing $f\in \mathcal{F}C_b^1$. By defining $\|d\rho\|$ and $\sigma_\rho(z)$ in the stated way, the extended relation (3.1) is exactly (3.3).$\hfill\Box$
\vskip.10in
By [Pu98, Theorem 4.1] (we will recall this result as Theorem A.1 in the Appendix) we can find a large class of sets $U\subset E$ such that $I_U\in BV(H,H_1)$ for $H_1$ properly chosen. This is contained in the following theorem. As preparation let us take $(Q^{1/2}H,\langle Q^{-1/2}\cdot,Q^{-1/2}\cdot\rangle)$ as the $(\bar{H},\langle \cdot,\cdot\rangle_{\bar{H}})$ in the Appendix. In this case $j_{\bar{H}}z=Qz$ for any $z\in E^*$ and $\nabla f(z)=QDf(z)$ for $f\in \mathcal{F}C_b^1$. Let $f$ satisfy the conditions in the Appendix, i.e. there exists a function $f\in H^{2,12}(E)$ such that $|Q^{-1/2}\nabla f|^{-1} \in L^{12}(E, \mu),  \nabla f \in Dom(\nabla^*)$, and set
 $$U:=f^{-1}((-\infty, 0)).$$
 We note that by the closability of $\nabla$, $\nabla f(z)\in \bar{H}=Q^{1/2}H$, hence $\nabla f(z)\in Dom(Q^{-1/2})$, where $Dom(Q^{-1/2})$ means the domain of $Q^{-1/2}$.
 Here  $H^{2,12}(E)$ is the completion of the space $$\mathcal{F}C_b^\infty:=\{u:u(z)=f({ }_{E^*}\!\langle l_1,z\rangle_E,{ }_{E^*}\!\langle l_2,z\rangle_E,...,{ }_{E^*}\!\langle l_m,z\rangle_E),z\in E, l_1,l_2,...,l_m\in E^*, m\in \mathbb{N}, f\in C_b^\infty(\mathbb{R}^m)\}.$$ with respect to the norm $$\|\varphi\|^{12}_{2,12}=\int(\varphi^{2}(x)+\sum_k\lambda_k(\frac{\partial\varphi}{\partial e_k})^2+\sum_{k,h}\lambda_k\lambda_h(\frac{\partial}{\partial e_k}\frac{\partial\varphi}{\partial e_h})^2)^6\mu(dx).$$

Set $\Sigma:=f^{-1}(0)$ and let $\nu$ be the corresponding surface measure constructed in [Pu98, Section 3]. Then by Theorem A.1, we have the following theorem.

 \th{Theorem 3.3} Assume there exists a function $f\in H^{2,12}(E)$ such that $|Q^{-1/2}\nabla f|^{-1} \in L^{12}(E, \mu),  \nabla f\in Dom(\nabla^*)$ and $\|Q^{-1}\nabla f(z)\|_{H^*_1}\nu(dz)$ is finite on $\Sigma$ for some separable Hilbert space $H_1$ such that $Q^{1/2}(H)\subset H_1\subset H$ (e.g. take $H_1:=Q^{1/2}(H)$). Then $I_U\in BV(H,H_1)$  and
 $$ \int_UD^*G(z)\mu(dz)=-\int_\Sigma { }_{H_1}\!\langle G(z),n_U(z)\rangle_{H_1^*}\|d\rho\|(dz)\textrm{ }\forall G\in (\mathcal{F}C_b^1)_{Q^{1/2}(H)},$$
 where $n_U(z)=Q^{-1}\nabla f(z)/\|Q^{-1}\nabla f(z)\|_{H^*_1}$ and $\|d\rho\|(dz)=\|Q^{-1}\nabla f(z)\|_{H^*_1}\nu(dz)$ is a finite measure on $\Sigma$.
Here we use $H$ also as a pivot space for $Q^{1/2}(H)^*$, i.e. $Q^{1/2}(H)\subset H_1\subset H(\equiv H^*) \subset H_1^*\subset (Q^{1/2}(H))^*$ and $Q^{-1}$ is considered as a continuous linear operator from $Q^{1/2}(H)$ to $(Q^{1/2}(H))^*$ (defined in the proof below).

\proof By Hypothesis 2.2 and  Theorem A.1 we have that  for any $G\in (\mathcal{F}C_b^1)_{Q^{1/2}(H)}$ $$\aligned \int_U\nabla^*G(z)\mu(dz)&=-\int_\Sigma \langle G(z),n(z)\rangle_{Q^{1/2}H}\mu_\sigma(dz)\\&=-\int_\Sigma \langle Q^{-1/2}G(z),Q^{-1/2}n(z)\rangle\mu_\sigma(dz)\\&=-\int_\Sigma {}_{Q^{1/2}(H)}\!\langle G(z),Q^{-1}n(z)\rangle_{(Q^{1/2}(H))^*}\mu_\sigma(dz)\\&=-\int_\Sigma {}_{Q^{1/2}(H)}\!\langle G(z),Q^{-1}\nabla f(z)\rangle_{(Q^{1/2}(H))^*}\nu(dz)\\&=-\int_\Sigma {}_{H_1}\!\langle G(z),Q^{-1}\nabla f(z)\rangle_{H_1^*}\nu(dz)\\&=-\int_\Sigma { }_{H_1}\!\langle G(z),n_U(z)\rangle_{H_1^*}\|d\rho\|(dz),\endaligned$$
 where $n(z)=\nabla f(z)/|Q^{-1/2}\nabla f(z)|$ and $\mu_\sigma(dz)=|Q^{-1/2}\nabla f(z)|\nu(dz)$ is a finite measure on $\Sigma$ and in the fifth equality we use that $\|Q^{-1}\nabla f(z)\|_{H^*_1}\nu(dz)$ is finite on $\Sigma$ and in the last step we take $n_U(z)$ and $\|d\rho\|$ as above.
 Here in the third equality we use $Q^{-1}$ as an operator from $Q^{1/2}(H)$ to $(Q^{1/2}(H))^*$. Indeed for any $y\in Q(H)$ we have
 $$\langle Q^{-1}y,z\rangle\leq |Q^{-1/2}y||Q^{-1/2}z|, \forall z\in Q^{1/2}(H),$$which implies that for any $y\in Q(H)$, we can define $\Lambda_y:=\langle Q^{-1}y,\cdot\rangle\in (Q^{1/2}(H))^*$. Clearly, $y\mapsto \Lambda_y$ is a continuous linear operator from $Q^{1/2}(H)$ to $(Q^{1/2}(H))^*$ with dense domain $Q(H)$. Hence it extends by continuity to all of $Q^{1/2}(H)$.
 By the definition of $D^*$, $\nabla^*$
we also have that for any $G\in (\mathcal{F}C_b^1)_{Q^{1/2}(H)}$, $D^*G=\nabla^*G$. Then by Theorem 3.2 the result follows.$\hfill\Box$
\vskip.10in
\th{Remark} From the definition of $(D^*,Dom(D^*))$ and $(\nabla^*,Dom(\nabla^*))$, it is easy to see that if we take $\bar{H}=Q^{1/2}(H)$ then $Dom(\nabla^*)\subset Dom(D^*)$ and $\nabla^*f=D^*f$ for $f\in Dom(\nabla^*)$. But if $f\in H^{1,2}(E)$, then $\nabla f\in Dom(D^*)$ implies that $\nabla f\in Dom(\nabla^*)$ since $\nabla f\in L^2(E,\mu;\bar{H})$. Here  $H^{1,2}(E)$ is the completion of the space $\mathcal{F}C_b^\infty$ with respect to the norm $$\|\varphi\|^{2}_{1,2}=\int\varphi^{2}(x)+\sum_k\lambda_k(\frac{\partial\varphi}{\partial e_k})^2\mu(dx).$$

\vskip.10in
 \th{Theorem 3.4} Let $\rho\in QR(E)\cap BV(H,H_1)$ and consider the measure $\|d\rho\|$ and $\sigma_\rho$ from Theorem 3.2(ii). Then for any smooth measure $\gamma$ under $P_\gamma:=\int P_zd\gamma$ there exists an $\mathcal{M}_t$- cylindrical Wiener process $W$, such that the sample paths of the associated distorted process $M^\rho$ on $F$ satisfy the following:  for $l\in  H_1\cap E^*\cap Q^{1/2}(H)$
 $${ }_{E^*}\!\langle l,X_t-X_0\rangle_{E}=\int_0^t\langle l,dW_s^z\rangle+\frac{1}{2}\int_0^t  { }_{H_1}\!\langle l,\sigma_\rho(X_s)\rangle_{H_1^*}dL_s^{\|d\rho\|}+\frac{1}{2}\int_0^t\beta_l(X_s) ds\textrm{ } \forall t\geq 0  \textrm{ }P_\gamma\rm{-a.s.}.\eqno(3.17)$$
 Here $L_t^{\|d\rho\|}$ is the real valued PCAF associated with $\|d\rho\|$ by the Revuz correspondence.

\proof Let $\{e_j\}$ be the orthonormal basis of H from Hypothesis 2.2. Define for all $k\in \mathbb{N}$ $$W_k(t):={ }_{E^*}\!\langle e_k,X_t-z\rangle_E-\frac{1}{2}\int_0^t  { }_{H_1}\!\langle e_k,\sigma_\rho(X_s)\rangle_{H_1^*}dL_s^{\|d\rho\|}-\frac{1}{2}\int_0^t\beta_{e_k}(X_s)ds .\eqno(3.18)$$ By (2.1) and (3.16) we  get for all $k\in \mathbb{N}$
 $$\mathcal{E}^\rho(e_k(\cdot),g)=-\frac{1}{2}\int_Eg(z)\beta_{e_k}(z)\rho(z)\mu(dz)-\frac{1}{2}\int_Eg(z)_{H_1}\langle e_k,\sigma_\rho(z)\rangle_{H_1^*}\|d\rho\|(dz)\textrm{ } \forall g\in \mathcal{F}C_b^1.$$
  By Theorem 2.4 it follows that for all $k\in \mathbb{N}$
 $$N_t^{e_k}=\frac{1}{2}\int_0^t { }_{H_1}\!\langle e_k,\sigma_\rho(X_s)\rangle_{H_1^*} dL_s^{\|d\rho\|}+\frac{1}{2}\int_0^t \beta_{e_k}(X_s) ds .\eqno(3.19)$$
 Here we get from (3.18), (3.19) and the uniqueness of decomposition (2.2) that,
 $$W_k(t)=M_t^{e_k}\textrm{ } \forall t\geq 0 \textrm{ }P_\gamma\rm{-a.s.} . $$
 By Dirichlet form theory we get $\langle M^{e_i},M^{e_j}\rangle_t=t\delta_{ij}$. So $W_k$ is an $\mathcal{M}_t$-Wiener process under $P_\gamma$. Thus, with $W$ being an $\mathcal{M}_t$- cylindrical Wiener process given by $W(t)=(W_k(t)e_k)_{k\in \mathbb{N}}$, (3.17) is satisfied $P_\gamma-a.e.$.    $\hfill\Box$
 \vskip.10in
\section{Examples }
BV functions in abstract Wiener space (see [Fu00], [FH01]) and BV functions in a Gelfand triple (see [RZZ12]) are examples of the extended notion of BV functions defined above, if the functions belong to $L^2(E,\mu)$. If we take $Q=I,$ and $\mu$ is a Gaussian measure with Cameron-Martin space $H=H_1=H_1^*$, they are just the BV functions in abstract Wiener space. If we take $H=E=E^*, Q=A^{-1}$ and $\mu$ is the Gaussian measure with mean zero and covariance operator $Q$, they are just the BV functions in a Gelfand triple. Now we want to give  examples which cannot be covered by the previous framworks, but only by the approach in the present paper.

\subsection{ The  stochastic reflected quantization equations in finite volume}

In this section we apply our BV functions theory to the stochastic quantization of $(\mathcal{P}(\phi)_2-)$ field theory in finite volume.
We consider the reflected problem in this case.
 Let $H=L^2(\Lambda;dx),$ where $\Lambda$ is a bounded rectangle in $\mathbb{R}^2$. Let $(-\Delta + I)_N$, be the generator of the
following quadratic form on $L^2(\Lambda;dx) : (u,v)\mapsto \int_\Lambda \langle\nabla u, \nabla v\rangle_{\mathbb{R}^d} dx + \int_\Lambda uv dx$
with $u, v\in \{g\in L^2(\Lambda;dx)|\nabla g\in L^2(\Lambda; dx)\}$ (where $\nabla$ is in the sense of
distributions). Let $\{e_n| n\in\mathbb{N}\}\subset C^\infty(\bar{\Lambda})$ be the (orthonormal) eigenbasis of
$(-\Delta+I)_N$ in $H$ and $\{\alpha_n|n\in\mathbb{N}\}\subset (0,\infty)$ the corresponding eigenvalues. Define for $\alpha\in\mathbb{R}$,
$$H^\alpha:=\{u\in L^2(\Lambda,dx) |\sum_{n=1}^\infty\alpha_n^\alpha\langle u,e_n\rangle^2_{L^2(\Lambda,dx)}<\infty\},$$
equipped with the inner product
$$\langle u,v\rangle_{H^\alpha}:=\sum_{n=1}^\infty\alpha_n^\alpha\langle u,e_n\rangle_{L^2(\Lambda,dx)}\langle v,e_n\rangle_{L^2(\Lambda,dx)}.$$
Set $E=H^{-s}, E^*=H^s$ for some $s>0$. Also set $\mu_0=N(0,(-\Delta+1)_N^{-1}):=N(0,C)$. Then $\mu_0$ is a measure supported on $E$.

For $h\in H^{-1}$ we define $X_h\in L^2(E,\mu_0)$ by $X_h:=\lim_{n\rightarrow\infty}{ }_{E^*}\!\langle k_n,\cdot\rangle_{E}\textrm{ in } L^2(E,\mu_0)$ where $k_n$ is any sequence in $E^*$ such that $k_n\rightarrow h$ in $H^{-1}$. We have the well-known (Wiener-It\^{o}) chaos decomposition $$L^2(E,\mu_0)=\bigoplus_{n\geq0}\mathcal{H}_n.$$
For $h\in L^2(\Lambda,dx)$ and $n\in \mathbb{N}$, define $:z^n:(h)$ to be the unique element in $\mathcal{H}_n$ such that
$$\int:z^n:(h):\prod_{j=1}^nX_{k_j}:_nd\nu=n!\int_{\mathbb{R}^2}\prod_{j=1}^n(\int_{\mathbb{R}^2}(-\triangle+1)_N^{-1}(x-y_j)k_j(y_j)dy_j)h(x)dx$$
where $k_1,...,k_n\in E^*$ and $:$ $:_n$  means orthogonal projection onto $\mathcal{H}_n$ (see [S74, V.1] for existence of $:z^n:(h)$). By [R86, Theorem 3.1] $h\in H\rightarrow :z^n:(h)\in L^2(E,\mu_0)$ is continuous. So by [AR91, Proposition 6.9], there exists a $\mathcal{B}(H^{-2})/ \mathcal{B}(H^{-2})$ measurable map $:z^n::H^{-2}\rightarrow H^{-2}$ such that $:z^n:(h)=_{H^{-2}}\langle :z^n:,h\rangle_{H^{2}}$. Finally, we set $:P(z):=\sum_{n=0}^{2N}a_n:z^n:$. Now we assume that $a_n\in \mathbb{R}$ and $a_{2N}>0$.

Let $$\mu=\exp{(-\int_{\Lambda}:P(z):dz)}\mu_0.$$

Then by [RZ92, Section 7], $\varphi(z):=\exp{(-\int_{\Lambda}:P(z):dz)}\in L^p(E,\mu_0), p\in[1,\infty)$.
Now define $Q:H\rightarrow H$ be a symmetric linear operator satisfying $Qe_k:=\frac{1}{\alpha_k^{2+s}}e_k, k\in \mathbb{N}$. Then $Q^{1/2}(H)=H^{2+s}$. Thus by [GlJ86, (9.1.32)] we have the following:

\th{Theorem 4.1.1} $Q^{1/2}(H)\subset H(\mu)$. Moreover for each $l\in Q^{1/2}H$, we have$$\beta_l(z)=_{H^{-s-2}}\!\langle -\sum_{n=1}^{2N}na_n:z^{n-1}:,l\rangle_{H^{2+s}}+_{H^{s}}\!\langle \Delta l-l,z\rangle_{H^{-s}}.$$
 \vskip.10in
Thus $Q$ satisfies Hypothesis 2.2. Now fix $k\in \mathbb{N}, a\in \mathbb{R}$ and take $U=\{z\in H^{-s}:{}_{H^{-s}}\!\langle z,e_k\rangle_{H^s}\leq a\},\rho=I_U$. Now take $f(z)={}_{H^{-s}}\!\langle z,e_k\rangle_{H^s}- a$. Then  $f\in L^p(E,\mu)$ for all $p\in [1,\infty)$. It is easy to get that $\nabla f(z)=Qe_k$ and all the conditions in Theorem 3.3 are satisfied.
Then by Theorem 3.3, $\rho$ is a BV function with $H=H_1=H^*_1$. Since $U$ is a convex closed set, it follows by [AR90, Theorem 3.2] that $\rho\in QR(E)$. Thus we can apply Theorem 3.4 directly and get the following:

\th{Theorem 4.1.2} There is an  $\mathcal{E}^\rho$-exceptional set $S\subset F$ such that $\forall z\in F\backslash S$ under $P_z$ there exists an $\mathcal{M}_t$- cylindrical Wiener process $W^z$, such that the sample paths of the associated distorted process $M^\rho$ on $F$ satisfy the following:  for $l\in H^{2+s}$
 $$\aligned{ }_{E^*}\!\langle l,X_t-X_0\rangle_{E}&=\int_0^t\langle l,dW_r^z\rangle-\frac{1}{2}\int_0^t \langle l,n_U(X_r)\rangle dL_r^{\|d\rho\|}\\+&\frac{1}{2}\int_0^t{ }_{H^{-s-2}}\!\langle -\sum_{n=1}^{2N}na_n:X_r^{n-1}:,l\rangle_{H^{2+s}}+_{H^{s}}\!\langle \Delta l-l,X_r\rangle_{H^{-s}} dr\textrm{ } \forall t\geq 0  \textrm{ }P_z\rm{-a.s.}.\endaligned$$
 Here $L_t^{\|d\rho\|}$ is the real valued PCAF associated with $\|d\rho\|$ by the Revuz correspondence, and $$I_{\partial U}(X_r)dL_r^{\|d\rho\|}=dL_r^{\|d\rho\|} \textrm{ } P-a.s.,$$ $n_U(z)=e_k$ is the normal to $\Sigma$.
 \vskip.10in
We can also take $U=\{z\in E:\|z\|^2_E\leq1\},\rho=I_U$. Now take $f(z)=\|z\|^2_E-1$. Then  $f\in L^p(E,\mu)$ for all $p\in [1,\infty)$. It is easy to get that $\nabla f(z)=2(-\Delta+1)_N^{-2s-2}z$. Thus we have $$\int|Q^{-1}\nabla f(z)|^{-12}\mu(dz)\leq (\int|Q^{-1}\nabla f(z)|^{-24}\mu_0(dz))^{1/2}(\int\varphi(z)^2\mu_0(dz))^{1/2}<\infty.$$
Since $|Q^{-1}\nabla f(z)|^{-1}\leq c(\sum_{k=1}^{25}{}_{E}\!\langle z,e_k\rangle_{E^*}^2)^{-1/2}$ with $c=1/\min\{2\alpha_k^{-s},k=1,...,25\}$,
$$\int|Q^{-1}\nabla f(z)|^{-24}\mu_0(dz)\leq c\int(\sum_{k=1}^{25}{}_{E}\!\langle z,e_k\rangle_{E^*}^2)^{-12}\mu_0(dz)<\infty,$$because the last integral transforms into an integral with respect to a Gaussian measure on $\mathbb{R}^{25}$ with an integrand having a singularity at zero of type $|x|^{-24}$, which is thus integrable.
 Also by $|Q^{-1}\nabla f(z)|\leq 2\|z\|_E$ we get $|Q^{-1}\nabla f(z)|\nu(dz)$ is finite on $\Sigma$.
 Hence
 all the conditions in Theorem 3.3 are satisfied.
Then by Theorem 3.3, $\rho$ is a BV function with $H_1=H_1^*=H$. Since $U$ is a convex closed set, as above we have $\rho\in QR(E)$. Thus as Theorem 4.1.2, we get the following:

\th{Theorem 4.1.3} There is an  $\mathcal{E}^\rho$-exceptional set $S\subset F$ such that $\forall z\in F\backslash S$ under $P_z$ there exists an $\mathcal{M}_t$- cylindrical Wiener process $W^z$, such that the sample paths of the associated distorted process $M^\rho$ on $F$ satisfy the following:  for $l\in H^{2+s}$
 $$\aligned{ }_{E^*}\!\langle l,X_t-X_0\rangle_{E}&=\int_0^t\langle l,dW_r^z\rangle-\frac{1}{2}\int_0^t \langle l,n_U(X_r)\rangle  dL_r^{\|d\rho\|}\\+&\frac{1}{2}\int_0^t{ }_{H^{-2-s}}\!\langle -\sum_{n=1}^{2N}na_n:X_r^{n-1}:,l\rangle_{H^{2+s}}+_{H^{s}}\!\langle \Delta l-l,X_r\rangle_{H^{-s}} dr\textrm{ } \forall t\geq 0  \textrm{ }P_z\rm{-a.s.}.\endaligned$$
 Here $L_t^{\|d\rho\|}$ is the real valued PCAF associated with $\|d\rho\|$ by the Revuz correspondence, and $$I_{\partial U}(X_r)dL_r^{\|d\rho\|}=dL_r^{\|d\rho\|} \textrm{ } P-a.s.,$$ $n_U(x)=\frac{(-\Delta+1)_N^{-s}x}{|(-\Delta+1)_N^{-s}x|}$ .
 \vskip.10in
Now we want to construct an example which is a BV functions in a Gelfand triple with $H\neq H_1$. Set $z_n=\sum_{k=1}^n\alpha_k^{-s/2}e_k$. Then it is obvious that $\langle z_n,\cdot\rangle$ converges to some function in $H^{2,12}(E,\mu)$. By the proof of [MR92, Proposition III.3.5] we can choose a subsequence $n_i, i\in \mathbb{N}$ such that $\langle z_{n_i},\cdot\rangle, i\in \mathbb{N}$ converges $Cap_{1,12}$- quasi-uniformly. Defining $z(x)=\limsup_{i\rightarrow\infty}\langle z_{n_i},x\rangle, x\in E,$ we obtain a $Cap_{1,12}$-quasi-continuous version of this function. We take $U=\{x\in E:z(x)\leq a\}$ for some $a\in \mathbb{R}$ such that $\mu(U)>0$, and $\rho=I_U$. Now take $f(x)=z(x)-a$. It is easy to get that $\nabla f(x)=\sum_k\alpha_k^{-2-\frac{3s}{2}}e_k$ and all the conditions in Theorem 3.3 are satisfied. Then by Theorem 3.3, $\rho$ is a BV function with $H_1=H^1, H_1^*=H^{-1}$. Since $z(x+re_k)=z(x)+r\alpha_k^{-s/2}$, it is continuous in $r$, by [AR90, Theorem 3.2] we have $\rho\in QR(E)$. Thus as Theorem 4.1.2, we get the following:

\th{Theorem 4.1.4} There is an  $\mathcal{E}^\rho$-exceptional set $S\subset F$ such that $\forall z\in F\backslash S$ under $P_z$ there exists an $\mathcal{M}_t$- cylindrical Wiener process $W^z$, such that the sample paths of the associated distorted process $M^\rho$ on $F$ satisfy the following:  for $l\in H^{2+s}$
 $$\aligned{ }_{E^*}\!\langle l,X_t-X_0\rangle_{E}&=\int_0^t\langle l,dW_r^z\rangle-\frac{1}{2}\int_0^t { }_{H^1}\!\langle l,n_U(X_r)\rangle_{H^{-1}} dL_r^{\|d\rho\|}\\+&\frac{1}{2}\int_0^t{ }_{H^{-2-s}}\!\langle -\sum_{n=1}^{2N}na_n:X_r^{n-1}:,l\rangle_{H^{2+s}}+_{H^{s}}\!\langle \Delta l-l,X_r\rangle_{H^{-s}} dr\textrm{ } \forall t\geq 0  \textrm{ }P_z\rm{-a.s.}.\endaligned$$
 Here $L_t^{\|d\rho\|}$ is the real valued PCAF associated with $\|d\rho\|$ by the Revuz correspondence, and $$I_{\partial U}(X_r)dL_r^{\|d\rho\|}=dL_r^{\|d\rho\|} \textrm{ } P-a.s.,$$ $n_U(x)=\frac{\sum_k\alpha_k^{-s/2}e_k}{\|\sum_k\alpha_k^{-s/2}e_k\|_{H^{-1}}}$.
 \vskip.10in
\th{Remark 4.1.5}From the above three theorems, we get a martingale solution to the stochastic reflected quantization equations. Choosing $a_n=0$ for all $0\leq n\leq 2N$, as a special case obtain the martingale solution to the stochastic reflected OU equations in two space dimensions.

\subsection{The reflected stochastic  quantization equations in infinite volume}

In this section, we consider the  reflected stochastic quantization equations in infinite volume. Let $\mathcal{S}'(\mathbb{R}^2)$ be the space of tempered Schwartz distributions on $\mathbb{R}^2$ and $\mathcal{S}(\mathbb{R}^2)$ the associated test function space equipped with the usual topology. Let $\mu_0$ be the mean zero Gaussian measure on $(\mathcal{S}'(\mathbb{R}^2),\mathcal{B}(\mathcal{S}'(\mathbb{R}^2)))$ with covariance
$$\int {}_{\mathcal{S}}\!\langle k_1,z\rangle_{\mathcal{S}'}{}_{\mathcal{S}}\!\langle k_2,z\rangle_{\mathcal{S}'}\mu_0(dz)=\int\int(-\Delta+1)^{-1}(x-y)k_1(x)k_2(y)dxdy=:\langle k_1,k_2\rangle_{H^{-1}},$$ where $(-\Delta+1)^{-1})$ denotes the Green function of the operator $(-\Delta+1)$ on $\mathbb{R}^2$. Let $H^{-1}$ be the real Hilbert space obtained
by completing $S(\mathbb{R}^2)$ w.r.t, the norm associated with the inner product $\langle \cdot,\cdot\rangle_{H^{-1}}$. Now for $n\in \mathbb{N}$, let $\mathcal{S}_{-n}$ denote the Hilbert subspace of $\mathcal{S}'(\mathbb{R}^2)$ which is the dual of $\mathcal{S}_n$ defined as the completion of $\mathcal{S}(\mathbb{R}^2)$ w.r.t the norm
$$\|k\|_n:=[\sum_{|m|\leq n}\int_{\mathbb{R}^2}(1+|x|^2)^n|(\frac{\partial^{m_1}}{\partial x_1^{m_1}},\frac{\partial^{m_2}}{\partial x_2^{m_2}})k(x)|^2dx]^{1/2}.$$For $h\in H^{-1}$ we define $X_h\in L^2(\mathcal{S}'(\mathbb{R}^2),\mu_0)$ by $X_h:=\lim_{n\rightarrow\infty}{ }_{\mathcal{S}}\!\langle k_n,\cdot\rangle_{\mathcal{S}'}\textrm{ in } L^2(\mathcal{S}'(\mathbb{R}^2),\mu_0)$ where $k_n$ is any sequence in $\mathcal{S}(\mathbb{R}^2)$ such that $k_n\rightarrow h$ in $H^{-1}$. We have the well-known (Wiener-It\^{o}) chaos decomposition $$L^2(\mathcal{S}'(\mathbb{R}^2),\mu_0)=\bigoplus_{n\geq0}\mathcal{H}_n.$$
For $h\in L^2(\mathbb{R}^2,dx)$ and $n\in \mathbb{N}$, define $:z^n:(h)$ to be the unique element in $\mathcal{H}_n$ such that
$$\int:z^n:(h):\prod_{j=1}^nX_{k_j}:_nd\nu=n!\int_{\mathbb{R}^2}\prod_{j=1}^n(\int_{\mathbb{R}^2}(-\Delta+1)^{-1}(x-y_j)k_j(y_j)dy_j)h(x)dx$$
where $k_1,...,k_n\in \mathcal{S}(\mathbb{R}^2)$ and $:$ $:_n$  means orthogonal projection onto $\mathcal{H}_n$ (see [S74, V.1] for existence of $:z^n:(h)$).

From now on we fix $N\in \mathbb{N}, a_n\in \mathbb{R}, 0\leq n\leq 2N$, and define for $h\in L^2(\mathbb{R}^2,dx)$
$$:P(z):(h):=\sum_{n=0}^{2N}a_n:z^n:(h)\textrm{ with } a_{2N}>0.$$We have that $\exp(-:P(z):(h))\in L^p(\mathcal{S}'(\mathbb{R}^2),\nu)$ for all $p\in [1,\infty)$ if $h\geq0$ (cf. [AR91, Section 7]), hence the following probability measures (called space-time cut-off quantum fields) are well-defined for $\Lambda\in \mathcal{B}(\mathbb{R}^2), \Lambda$ bounded,
$$\mu_\Lambda:=\frac{\exp{(-:P(z):(1_\Lambda))}}{\int\exp{(-:P(z):(1_\Lambda))}d\mu_0}\mu_0.$$
It has been proven that the weak limit $$\lim_{\Lambda\rightarrow \mathbb{R}^2}\mu_\Lambda=:\mu$$ exists as a probability measure on $(\mathcal{S}'(\mathbb{R}^2),\mathcal{B}(\mathcal{S}'(\mathbb{R}^2)))$ having moments of all orders (see [GlJ86] and also [AR91, Section 7]). In particular, it follows by [AR89, Proposition 3.7] that $\mu(\mathcal{S}_{-n})=1$ for $n\in \mathbb{N}$ large enough. Now we take $E:=\mathcal{S}_{-n}, H:=L^2(\mathbb{R}^2,dx)$. By [AR91, Theorem 7.11] we know that $C_0^\infty(\mathbb{R}^2)\subset H(\mu)$ and that for all $k\in C_0^\infty(\mathbb{R}^2)$ we have
$$\beta_k(z):=-\sum_{m=1}^{2N}ma_m:z^{m-1}:(k)-{}_{\mathcal{S}_n}\!\langle (-\Delta+1)k,z\rangle_{\mathcal{S}_{-n}},z\in \mathcal{S}_{-n},$$ and that, choosing $n$ large, if necessary, there exists $\beta:\mathcal{S}_{-n}\rightarrow \mathcal{S}_{-n}$, $\mathcal{B}(\mathcal{S}_{-n})/\mathcal{B}(\mathcal{S}_{-n})$-measurable, such that for all $k\in C_0^\infty(\mathbb{R}^2)$
$$\beta_k={}_{\mathcal{S}_n}\!\langle k,\beta\rangle_{\mathcal{S}_{-n}}\qquad \mu-a.e.$$
and $$\int_{\mathcal{S}_{-n}}\|\beta\|^2_{\mathcal{S}_{-n}}<\infty.$$
Defining $\beta_k:=\langle k,\beta\rangle_{S_{-n}}, k\in \mathcal{S}_{-n}$, because $C_0^\infty(\mathbb{R}^2)$ is dense in $\mathcal{S}_n$, it follows by Lebesgue's dominated convergence theorem that $\mu$ on $E:=\mathcal{S}_{-n}$ is differentiable along each $k\in \mathcal{S}_n$ with $\beta_k\in L^2(E,\mu)$, i.e. $\mathcal{S}_n\subset H(\mu)$. By e.g [H80, A.3], the embedding $H\subset E$ is Hilbert-Schmidt. Let us denote the inclusion map $H\subset E$ by $J$ and let $J^*:E\rightarrow H$ be its adjoint. Define $Q_1:=JJ^*$. Then $\ker Q_1=\{0\}$, since $H$ is dense in $E$, and by [PR07, Proposition 2.5.2] $Q_1$ is a nonngegative definite and symmetric bounded linear operator with finite trace. Furthermore, $J:H\rightarrow Q_1^{1/2}E$ is an isometry, i.e. $\forall h_1,h_2\in H$
$$\langle h_1,h_2\rangle=\langle Q^{-1/2}_1h_1,Q^{-1/2}_1h_2\rangle_E.$$
In particular, $Q:=Q_1|_H$ is symmetric and nonnegative definite on $H$ with finite trace. It is elementary to check that then $Q^{1/2}H=\mathcal{S}_n$. In particular, $Q$ has an extension, namely $Q_1$, from $E$ to $H$, and thus $Q$ is an operator satisfying Hypothesis 2.2, i.e.
 there exists an orthonormal basis $e_k,k\in \mathbb{N},$ of $H$ and $\lambda_k\in (0,\infty)$ such that $Qe_k=\lambda_ke_k$ and $\frac{1}{\sqrt{\lambda_k}}e_k,k\in \mathbb{N},$ is an orthonormal basis of $E$. Then by the same arguments as in the last section, we obtain the following two theorems.

Fix $k\in \mathbb{N}, a\in \mathbb{R}$ and take $U=\{z\in \mathcal{S}_{-n}: {}_{\mathcal{S}_{-n}}\!\langle z,e_k\rangle_{\mathcal{S}_{n}}\leq a\},\rho=I_U$. Now take $f(z)={}_{\mathcal{S}_{-n}}\!\langle z,e_k\rangle_{\mathcal{S}_{n}}- a$. Then  $f\in L^p(E,\mu)$ for all $p\in [1,\infty)$. It is easy to get that $\nabla f(z)=Qe_k$ and all the conditions in Theorem 3.3 are satisfied. Then by Theorem 3.3, $\rho$ is a BV function with $H=H_1=H^*_1$. Since $U$ is a convex closed set, as above we have $\rho\in QR(E)$. Thus we can apply Theorem 3.4 directly and  get the following:

\th{Theorem 4.2.1} There is an  $\mathcal{E}^\rho$-exceptional set $S\subset F$ such that $\forall z\in F\backslash S$ under $P_z$ there exists an $\mathcal{M}_t$- cylindrical Wiener process $W^z$, such that the sample paths of the associated distorted process $M^\rho$ on $F$ satisfy the following:  for $l\in \mathcal{S}(\mathbb{R}^2)$
 $$\aligned{ }_{E^*}\!\langle l,X_t-X_0\rangle_{E}&=\int_0^t\langle l,dW_s^z\rangle-\frac{1}{2}\int_0^t \langle l,n_U(X_s)\rangle dL_s^{\|d\rho\|}\\+&\frac{1}{2}\int_0^t -\sum_{n=1}^{2N}na_n:X_r^{n-1}:(l)+_{\mathcal{S}_{n}}\!\langle \Delta l-l,X_r\rangle_{\mathcal{S}_{-n}} dr\textrm{ } \forall t\geq 0  \textrm{ }P_z\rm{-a.s.}.\endaligned$$
 Here $L_t^{\|d\rho\|}$ is the real valued PCAF associated with $\|d\rho\|$ by the Revuz correspondence, and $$I_{\partial U}(X_s)dL_s^{\|d\rho\|}=dL_s^{\|d\rho\|} \textrm{ } P-a.s.,$$ $n_U(z)=e_k$ is the normal to $\Sigma$.
 \vskip.10in
We can also take $U=\{z\in E: \|z\|^2_E\leq1\},\rho=I_U$. Now take $f=\|z\|^2_E-1$, and then  $f\in L^p(E,\mu)$ for all $p\in [1,\infty)$. It is easy to get that $\nabla f(z)=2Q^2z\in Q^{1/2}(H)$. Then by [R86, Proposition 6.8] we have that
 $$\int(|Q^{-1}\nabla f(z)|^{-12}\wedge N)\mu_\Lambda(dz)\leq C\int|Q^{-1}\nabla f(z)|^{-12}\mu_0(dz)<\infty.$$
Taking limit in $N$ and $\Lambda$ we get that $$\int|Q^{-1}\nabla f(z)|^{-12}\mu(dz)\leq C\int|Q^{-1}\nabla f(z)|^{-12}\mu_0(dz)<\infty.$$ Also by $|Q^{-1}\nabla f(z)|\leq \|z\|_E$ we get $|Q^{-1}\nabla f(z)|\nu(dz)$ is finite on $\Sigma$. Hence all the conditions in Theorem 3.3 is satisfied and thus $\rho$ is an BV function with $H_1=H_1^*=H$. Since $U$ is a convex closed set, we have $\rho\in QR(E)$. Thus as Theorem 4.1.2, we get the following:

\th{Theorem 4.2.2} There is an  $\mathcal{E}^\rho$-exceptional set $S\subset F$ such that $\forall z\in F\backslash S$ under $P_z$ there exists an $\mathcal{M}_t$- cylindrical Wiener process $W^z$, such that the sample paths of the associated distorted process $M^\rho$ on $F$ satisfy the following:  for $l\in \mathcal{S}(\mathbb{R}^2)$
 $$\aligned{ }_{E^*}\!\langle l,X_t-X_0\rangle_{E}&=\int_0^t\langle l,dW_s^z\rangle-\frac{1}{2}\int_0^t { }_{E^*}\!\langle l,n_U(X_s)\rangle _{E} dL_s^{\|d\rho\|}\\+&\frac{1}{2}\int_0^t -\sum_{n=1}^{2N}na_n:X_r^{n-1}:(l)+_{\mathcal{S}_n}\!\langle \Delta l-l,X_r\rangle_{\mathcal{S}_{-n}} dr\textrm{ } \forall t\geq 0  \textrm{ }P_z\rm{-a.s.}.\endaligned$$
 Here $L_t^{\|d\rho\|}$ is the real valued PCAF associated with $\|d\rho\|$ by the Revuz correspondence, and $$I_{\partial U}(X_s)dL_s^{\|d\rho\|}=dL_s^{\|d\rho\|} \textrm{ } P-a.s.,$$ $n_U(x)=\frac{\sum_k\lambda_k\langle x,e_k\rangle e_k}{|\sum_k\lambda_k\langle x,e_k\rangle e_k|}$ is the normal to $\Sigma$.
 \vskip.10in

\subsection{ Other examples}
 Consider $\mu=\varphi^2\mu_0$ with $\varphi\in L^2(H,\mu_0)$. Here $\mu_0$ is a Gaussian measure on $H$ satisfying Hypothesis 2.1 in [RZZ12], i.e. $A:D(A)\subset H\rightarrow H$ is a linear self-adjoint operator on H such that $\langle Ax,x\rangle\geq\delta|x|^2 \textrm{ }\forall x\in D(A)$ for some $\delta>0$ and $A^{-1}$ is of trace class. Assume that $\varphi$ is Fr\'{e}chet differentiable such that for its Fr\'{e}chet derivative we have
 $$\int_H|D\varphi|^2d\mu_0<\infty, \varphi>0 .\eqno(4.1)$$ Since $\mu_0$ satisfies the log-Sobolev inequality, by Young's inequality we can deduce $$\varphi(x)\cdot\langle e_k,x\rangle\in L^2(H,\mu_0) ,\forall k\in \mathbb{N}.$$
 By [MR92, II.3.d], for $l\in D(A)$, we have $$\beta_l(z)=-2\langle Al,z\rangle+2\langle l,\frac{D\varphi(z)}{\varphi(z)}\rangle.$$

Then by Theorem 3.4 we get the following result.

 \th{Theorem 4.3.1}Let $\rho\in QR(H)\cap BV(H,H_1)$ and consider the measure $\|d\rho\|$ and $\sigma_\rho$ from Theorem 3.2(ii). Then there is an  $\mathcal{E}^\rho$-exceptional set $S\subset F$ such that $\forall z\in F\backslash S$ under $P_z$ there exists an $\mathcal{M}_t$- cylindrical Wiener process $W^z$, such that the sample paths of the associated distorted process $M^\rho$ on $F$ satisfy the following:  for $l\in D(A)\cap H_1$
 $$\langle l,X_t-X_0\rangle=\int_0^t\langle l,dW_s^z\rangle+\frac{1}{2}\int_0^t  { }_{H_1}\!\langle l,\sigma_\rho(X_s)\rangle_{H_1^*}dL_s^{\|d\rho\|}-\int_0^t\langle Al, X_s\rangle ds+\int_0^t\langle l, \frac{D\varphi(X_s)}{\varphi(X_s)}\rangle ds \textrm{ }P_z\rm{-a.s.}.$$
 Here $L_t^{\|d\rho\|}$ is the real valued PCAF associated with $\|d\rho\|$ by the Revuz correspondence.
 \vskip.10in
Assume $f$ satisfies the same conditions as in Theorem 3.3 and $$U=f^{-1}((-\infty, 0)).$$

\th{Theorem 4.3.2} Let $I_U\in QR(H)$  satisfy the conditions in Theorem 3.3 and let $|Df|$ be finite on $\partial U$. Then there is an  $\mathcal{E}^\rho$-exceptional set $S\subset F$ such that $\forall z\in F\backslash S$ under $P_z$ there exists an $\mathcal{M}_t$- cylindrical Wiener process $W^z$, such that the sample paths of the associated distorted process $M^\rho$ on $F$ satisfy the following:  for $l\in D(A)$
 $$\langle l,X_t-X_0\rangle=\int_0^t\langle l,dW_s^z\rangle-\frac{1}{2}\int_0^t  \langle l,n_U(X_s)\rangle dL_s^{\|d\rho\|}-\int_0^t\langle Al, X_s\rangle ds+\int_0^t\langle l, \frac{D\varphi(X_s)}{\varphi(X_s)}\rangle ds \textrm{ }P_z\rm{-a.s.}.$$
 Here $L_t^{\|d\rho\|}$ is the real valued PCAF associated with $\|d\rho\|$ by the Revuz correspondence,   $n_U(z)=Df(z)/|Df(z)|$ is the normal to $\Sigma$.
 \vskip.10in
Now consider the following stochastic differential inclusion in the Hilbert space $H$,
  $$\left\{\begin{array}{ll}dX(t)+(AX(t)-\frac{D\varphi(X_t)}{\varphi(X_t)}+N_U(X(t)))dt\ni dW(t),&\ \ \ \ \textrm{ }\\X(0)=x,&\ \ \ \ \textrm{ } \end{array}\right.\eqno(4.2)$$
  where  $W(t)$ is a cylindrical Wiener process in $H$ on a filtered probability space $(\Omega,\mathcal{F},\mathcal{F}_t,P)$ and $N_U(x)$ is the normal cone to $U$ at $x$, i.e.
  $$N_U(x)=\{z\in H:\langle z,y-x\rangle\leq0\textrm{ }\forall y\in U\}.$$
  \vskip.10in
\th{Definition 4.3.3} A pair of continuous $H\times \mathbb{R}$-valued and $\mathcal{F}_t$-adapted processes $(X(t),L(t)),t\in [0,T]$, is called a solution of (4.2) if the following conditions hold.

(i) $X(t)\in U$ for all $t \in [0,T]$ $P-a.s.$;

(ii) $L$ is an increasing process with the property that
$$I_{\partial U}(X_s)dL_s=dL_s \textrm{ } P-a.s.$$ and for any $l\in D(A)$ we have
$$\langle l,X_t-x\rangle=\int_0^t\langle l,dW_s\rangle-\int_0^t\langle l,\textbf{n}_U(X_s)dL_s\rangle-\int_0^t\langle Al, X_s\rangle ds+\int_0^t\langle l, \frac{D\varphi(X_s)}{\varphi(X_s)}\rangle ds\textrm{ }\forall t\geq0 \textrm{ }P-a.s.$$
where $\textbf{n}_U$ is the exterior normal to $U$.
\vskip.10in
We recall that if $\log\varphi$ is concave, then $\langle\frac{D\varphi(x)}{\varphi(x)}-\frac{D\varphi(y)}{\varphi(y)},x-y\rangle\leq0$, for all $x,y\in H$. Hence by a modification of [RZZ12, Theorem 5.11], we obtain  pathwise uniqueness.

\th{Theorem 4.3.4}  Assume $U\subset H$ satisfies the same conditions as in Theorem 4.3.2, and $\log \varphi$ is a concave function. Then the stochastic inclusion (4.2) admits at most one solution in the sense of Definition 4.3.3.

\vskip.10in
Combining Theorem 4.3.2 and 4.3.4 with the Yamada-Watanabe Theorem, we now obtain the following:
\vskip.10in
\th{Theorem 4.3.5} Assume $U\subset H$ satisfies the same conditions as in Theorem 4.3.2, and that $\log \varphi$ is a concave function. Then there exists a Borel set $M\subset H$ with $I_U\cdot\mu(M)=\mu(U)$ such that for every $x\in M$, (4.2) has a pathwise unique continuous strong solution in the sense that for every probability space $(\Omega,\mathcal{F},\mathcal{F}_t,P)$ with an $\mathcal{F}_t$-Wiener process $W$,  there exists a unique pair of $\mathcal{F}_t$-adapted processes $(X,L)$ satisfying Definition 4.3.3. Moreover $X(t)\in M$ for all $t\geq0 $ $P$-a.s.
\vskip.10in
As an example, we can take $f=\langle x,x\rangle-1,\varphi(x)=e^{-|x|^4}$. Then $\log \varphi$ is a concave function and we can check that all the conditions in Theorem 4.3.5 are satisfied. Hence by  Theorem 4.3.5, we get that there exists a unique probabilistically strong solution in the sense of Definition 4.3.3 for the following problem:

 $$\left\{\begin{array}{ll}dX(t)+(AX(t)+4|X_t|^2X_t+N_U(X(t)))dt\ni dW(t),&\ \ \ \ \textrm{ }\\X(0)=x.&\ \ \ \ \textrm{ } \end{array}\right.$$

\section{Appendix}

In this appendix we will recall the Gauss-Ostrogradski\u{l} formula proved in [Pu98].  Let $E$ be a  separable Banach space and $\bar{H}$ be a
Hilbert space ( with $\{f_k\}$ be an orthonormal basis) continuously and densely embedded in $E$. We denote by  $\langle\cdot,\cdot\rangle_{\bar{H}}$ the
scalar product in $\bar{H}$, and by $\|\cdot\|$ the norm in $\bar{H}$. The map $ j_{\bar{H}}:E^*\mapsto \bar{H}$ defined by the formula
$$\langle h,j_{\bar{H}}(l) \rangle_{\bar{H}} = l(h) \qquad \forall h\in \bar{H}, l\in E^*,$$
is a continuous embedding of $E^*$ in $\bar{H}$. We also introduce a family of $\bar{H}$-valued functions on $E$ by
$$(\mathcal{F}C_b^1)_{\bar{H}}:=\{G:G(z)=\sum_{j=1}^m g_j(z)l^j,z\in E, g_j\in \mathcal{F}C_b^1,l^j\in \bar{H}\}$$

By $\nabla u$ denote the $\bar{H}$-derivative of $u\in \mathcal{F}C_b^1$, i.e. $\nabla u(x): = j_{\bar{H}}(u'(x))$, where $u'(x)\in E^*$ is the Fr\^{e}chet derivative of $u$ at $x\in E$. Since the operator $\nabla:\mathcal{F}C_b^1\subset L^2(E,\mu)\mapsto L^2(E,\mu;\bar{H})$ is closable, we can uniquely extend $\nabla$ to all of $W^{1,2}(E)$.
 By $\nabla^*$ denote the adjoint operator of $\nabla:\mathcal{F}C_b^1\in L^2(E,\mu)\mapsto L^2(E,\mu;\bar{H})$. That is
$$Dom(\nabla^*):=\{G\in L^2(E,\mu;\bar{H})|\mathcal{F}C_b^1\ni u\mapsto \int_E\langle G,\nabla u\rangle_{\bar{H}} d\mu \textrm{ is continuous with respect to } L^2(E,\mu)\}.$$

Let $\mu$ be a differentiable measure on $E$ in the sense of Definition 2.1 such that $\bar{H}\subset H(\mu)$. As in [Pu98, Section 4] let $f\in W^{2,12}(E)$ such that $\|\nabla f\|^{-1} \in L^{12}(E, \mu),  \nabla f \in Dom(\nabla^*)$. Since $f\in W^{2,12}(E)\subset W^{1,12}(E)$, $f$ has a $Cap_{1,12}$-quasi-continuous version which is again denoted by $f$. Set
 $$U:=f^{-1}((-\infty, 0)).$$ Here  $W^{1,12}(E), W^{2,12}(E)$ are the completion of the space $$\mathcal{F}C_b^\infty:=\{u:u(z)=f({ }_{E^*}\!\langle l_1,z\rangle_E,{ }_{E^*}\!\langle l_2,z\rangle_E,...,{ }_{E^*}\!\langle l_m,z\rangle_E),z\in E, l_1,l_2,...,l_m\in E^*, m\in \mathbb{N}, f\in C_b^\infty(\mathbb{R}^m)\}$$ with respect to the norm $$\|\varphi\|^{12}_{1,12}=\int(\varphi^{2}(x)+\sum_k(\frac{\partial\varphi}{\partial f_k})^2)^6d\mu$$
 and $$\|\varphi\|^{12}_{2,12}=\int(\varphi^{2}(x)+\sum_k(\frac{\partial\varphi}{\partial f_k})^2+\sum_{k,h}(\frac{\partial}{\partial f_k}\frac{\partial\varphi}{\partial f_h})^2)^6\mu(dx),$$
 respectively.
$Cap_{1,12}$ is defined as follows:$$Cap_{1,12}(U)=\inf\{\|\varphi\|_{1,12}:\varphi\geq1\textrm{ }\mu-a.e \textrm{ on }U\} \quad \textrm{for an open set }U,$$
$$Cap_{1,12}(A)=\inf\{Cap_{1,12}(U):U\textrm{ is open}, U\supset A\} \quad \textrm{for an arbitary set }A.$$

Set $\Sigma:=f^{-1}(0)$ and let $\nu$ be the corresponding surface measure constructed in [Pu98, Section 3]. Now we can restate [Pu98 , Theorem 4.1] in the following form:

\th{Theorem A.1} Assume $\bar{H}\subset H(\mu)$ and let $f$ be as above. Then the following  Gauss-Ostrogradski\u{l} formula holds for the $U$ defined above:
 $$\int_U\nabla^*G(z)\mu(dz)=-\int_\Sigma \!\langle G(z),n(z)\rangle_{\bar{H}}\mu_\sigma(dz)\qquad\forall G\in (\mathcal{F}C_b^1)_{\bar{H}},\eqno(A.1)$$
 where $n(z)=\nabla f(z)/\|\nabla f(z)\|$ and $\mu_\sigma(dz)=\|\nabla f(z)\|\nu(dz)$ is a finite measure on $\Sigma$.

\th{Remark A.2} (i) The formulation here is different from [Pu98, Theorem 4.1]. By [Bo10] Definition 6.6.1 and the same argument as [Bo10, Theorem 8.10.1] and the definition of $\nabla^*$, we can easily see that for any $G\in (\mathcal{F}C_b^1)_{\bar{H}}$, the divergence of $G$ exists and equals $\nabla^* G$ in $L^2(E,\mu)$ which implies (A.1).

(ii) Theorem A.1 (as [Pu98, Theorem 4.1]) does not depend on the $Cap_{1,12}$-quasi-continuous version chosen for $f$, because $U$ and $\Sigma$ will only change by $Cap_{1,12}$-zero sets, which have measure zero with respect to both $\mu$ and $\nu$ (see [Pu98, Lemma 3.2 ]).
\vskip 1cm
\th{Acknowledgement.} The authors would like to thank the referee for valuable comments and suggestions that led to a much improved version of this paper.

\end{document}